\documentclass[11pt]{article}
\usepackage{graphicx}
\usepackage{amsfonts}
\usepackage{bbm}
\usepackage{mathrsfs}
 \usepackage{mathrsfs,amsmath,amssymb,amsthm}
\usepackage[dvips]{color}
 \usepackage{amssymb}
 \usepackage{amsbsy}
 \usepackage[dvips]{color}
\allowdisplaybreaks
 \theoremstyle{plain}
\theoremstyle{remark}  \newtheorem{remark}{\noindent\mbox{Remark}}
 \theoremstyle{plain}
 \theoremstyle{plain}
\theoremstyle{plain} \newtheorem{theorem}{\noindent\mbox{Theorem}}
 \theoremstyle{plain}
 \theoremstyle{plain}
\theoremstyle{definition}

 \def\bq{\begin{equation}}
 \def\eq{\end{equation}}
 \def\eqn{\end{eqnarray}}
 \def\bqn{\begin{eqnarray}}
 
 \def\qed{\hfill$\Box$\medskip}
 \def\rto{\rightarrow\infty}
 \def\z{\left}
 \def\y{\right}
 \def\no{\nonumber}

 \def\ka{{\kappa}}

\begin{document}
 \title{\textbf{Stationary distribution for  birth and death process with one-side bounded jumps}\footnote{Supported by National
Nature Science Foundation of China (Grant No. 11226199) and Natural Science Foundation of Anhui Educational Committee (Grant No. KJ2014A085)}}                  

\author{  Hua-Ming \uppercase{Wang} \\
     Department of Mathematics, Anhui Normal University, Wuhu 241003, China\\
    E-mail\,$:$ hmking@mail.ahnu.edu.cn}
\date{}
\maketitle%

\vspace{-.5cm}

\begin{center}
\begin{minipage}[c]{12cm}
\begin{center}\textbf{Abstract}\quad \end{center}

In this paper, we study a birth and death process $\{N_t\}_{t\ge0}$ on positive half lattice, which at each discontinuity jumps at most a distance $R\ge 1$ to the right or exactly a distance $1$ to the left. The transitional probabilities at each site are nonhomogeneous. Firstly, sufficient conditions for the recurrence and positive recurrence are presented. Then by the branching structure within random walk with one-side bounded jumps set up in Hong and Wang (2013), the explicit form of the stationary distribution of the process $\{N_t\}_{t\ge0}$ is formulated.

\mbox{}\textbf{Keywords:}\  birth and death process; stationary distribution; embedded process; branching structure.\\
\mbox{}\textbf{MSC 2010:}\
 60J27;  60J80
\end{minipage}
\end{center}

\section{Introduction}

The aim of this paper is to study the birth and death process with one-side bounded jumps. Fix $1\le R\in \mathbb Z$ and let $(\mu_i,\lambda_i^1,...,\lambda_i^R)_{i\ge 0},$ where $\mu_0=0,$ be a sequence of nonnegative $\mathbb R^{R+1}$-valued vectors.
Let $\{N_t\}_{t\ge0}$ be a continuous time Markov chain, which waits at a state $n$ an exponentially distributed time with parameter $\mu_n+\sum_{r=1}^R\lambda^{r}_n$ and then jumps to $n+i$ with probability ${\lambda_n^i}/(\mu_n+\sum_{r=1}^R\lambda^{r}_n),$ $i=1,...,R$ or to $n-1$ with probability ${\mu_n}/(\mu_n+\sum_{r=1}^R\lambda^{r}_n).$  We call the process $\{N_t\}_{t\ge 0}$ a {\it nonhomogeneous birth and death process with one-side bounded jumps}.

The above defined process $\{N_t\}$ is the continuous time analogue of the so-called (1,R) random walk studied in Hong and Zhou \cite{hz}, where the stationary distribution of the walk was given.

In this paper, we present some sufficient conditions for recurrence and positive recurrence.  For the positive recurrent case, we formulate the explicit stationary distribution of $\{N_t\}.$ The idea to formulated the stationary distribution is as follows. By looking only at the discontinuities of $\{N_t\}$ we get its embedded process $\{X_n\},$ whose stationary distribution $\pi_k, k\ge 0$ could be formulated by mean of the branching structure constructed in Hong and Wang \cite{hwa}. But $\{\pi_k/(\mu_k+\sum_{r=1}^R\lambda^r_k)\}_{k\ge0}$ defines a invariant measure for the process $\{N_t\}.$ In this way, the stationary distribution for $\{N_t\}$ could be formulated.
The following condition is always assumed to be satisfied throughout the paper.

\noindent{\bf(C)}  There are small $\ka>0$ and large $K>0$ such that for all $n\ge0,$ $\ka< \mu_n+\sum_{r=1}^R\lambda_n^r<K.$

\begin{remark}

Let $Q=(q_{ij})_{i,j\ge0}$ be a matrix with $$q_{ij}=\left\{\begin{array}{ll}
                                       \mu_i,&\text{ if } j=i-1; \\
                                       \lambda_{i}^r,&\text{ if } j=i+r,\ l=1,...,R;\\
                                       -\big(\mu_i+\sum_{r=1}^R\lambda_i^r), &\text{ if } j=i;\\
                                       0, &\text{ else.}
                                     \end{array}
\right.$$
Then it follows from (C) that the matrix $Q$ is a conservative Q-matrix bounded from above. Hence the process $\{N_t\}$ exists (See for example Anderson \cite{ads}, Proposition 2.9, Chapter 2.). Some weaker condition could imply the existence, see for example Wang \cite{w}.
\end{remark}

For $i\ge 1,$ let $a_i^k=\frac{\sum_{l=k}^R\lambda_i^l}{\mu_i},$ $k=1,...,R,$
and introduce matrices
\begin{equation}M_i=\left(
        \begin{array}{cccc}
          a_{i}^1 &\cdots  & a_i^{R-1}& a_i^R \\
          1 & \cdots & 0 & 0 \\
          \vdots& \ddots & \vdots &\vdots  \\
          0 & \cdots & 1 & 0 \\
        \end{array}
      \right).
\end{equation}
Throughout the paper, for $1\le i\le R,$ $\mathbf{e}_i$ denotes the row vector whose $i$th component equals to $1$ and all other components equal to $0.$ Write also $\mathbf 1=(\mathbf{e}_1+...+\mathbf{e}_R)^T.$
Let $$P_t(i,j)=P(N_t=j|N_0=i).$$ The theorem below is the main result of the paper.
\begin{theorem}\label{main} Suppose that condition (C) holds. Then the following statements hold.

\noindent(a)  If $\lim_{n\rto}\mathbf{e}_1M_1M_2\cdots M_n\mathbf{e}_1^T=0,$ then $\{N_t\}$ is recurrent.

\noindent(b) If $\sum_{n=1}^\infty\frac{1}{\mu_i} \mathbf{e}_1M_1M_2\cdots M_n\mathbf{e}_1^T<\infty,$ then $\{N_t\}$ is positive recurrent and the limits \begin{equation}\label{lima}\psi_0:=\lim_{t\rto}P_t(i,0)=\frac{(\sum_{r=1}^R\lambda_0^r)^{-1}}{{(\sum_{r=1}^R\lambda_0^r)^{-1}}+\sum_{n=1}^\infty\frac{1}{\mu_n}\mathbf{e}_1M_1M_2\cdots M_{n-1}\mathbf{e}_1^T}\end{equation}
  \begin{equation}\label{limb}\psi_k:=\lim_{t\rto}P_t(i,k)=\frac{\frac{1}{\mu_n}\mathbf{e}_1M_1M_2\cdots M_{k-1}\mathbf{e}_1^T}{{(\sum_{r=1}^R\lambda_0^r)^{-1}}+\sum_{n=1}^\infty\frac{1}{\mu_n}\mathbf{e}_1M_1M_2\cdots M_{n-1}\mathbf{e}_1^T}\end{equation}
define a stationary distribution for $\{N_t\}$ in the sense that for all $t>0,$
$$\psi_k=\sum_{n=0}^\infty \psi_nP_t(n,k),\ k\ge 0.$$
\end{theorem}
\begin{remark}
We give a remark on the recurrence and positive recurrence of continuous time Markov chain $\{N_t\}.$ Let $\tau=\inf\{t>0:N_t\neq N_0\}$ and define $\eta=\inf\{t>\tau:N_t=N_0\}.$ Then by definition, $\tau$ is the time when $\{N_t\}$ leaves the starting state and $\eta$ is the time when $\{N_t\}$ returns to the starting state after it leaves it.
If $P(\eta<\infty)=1,$ we say that $\{N_t\}$ is recurrent; if $E\eta<\infty,$ we say that $\{N_t\}$ is positive recurrent.
\end{remark}

\section{Proof}\label{pro}

Note that under condition (C), $\{N_t\}$ exists. Let $\tau_0=0$ and define recursively for $n\ge 1,$ $$\tau_n=\inf\{t\ge\tau_{n-1}:N_t\neq N_{\tau_{n-1}}\}$$ with convention being that $\inf\phi=\infty.$ Set for $n\ge0,$ $X_n=N_{\tau_n}.$ Then $\{X_n\}_{n\ge0}$ is a discrete time Markov chain on positive half lattice with transition probabilities
\begin{equation}\label{rij}
  r_{ij}=\left\{\begin{array}{cl}
                  \frac{\mu_i}{\mu_i+\sum_{l=1}^R\lambda_i^l}, & j=i-1 \\
                 \frac{\lambda_i^k}{\mu_i+\sum_{l=1}^R\lambda_i^l}, & j=i+k, k=1,..,R\\
                 0,& \text{else.}
                \end{array}
\right.
\end{equation}
$\{X_n\}$ is also known as the {\it embedded process} of $\{N_t\}.$

For the embedded process $\{X_n\}$ define $$T=\inf\{k>0:X_k=0\},$$ which is the time $\{X_n\}$ returns to $0.$

For the  process $\{N_t\}$ define $$\eta=\inf\{t>\tau_1:N_t=N_0\},$$ which is the time it returns to $0$ after it leaves $0.$

To prove part (a) of the theorem, it suffices to show that $\{X_n\}$ is recurrent whenever  $\lim_{n\rto}\mathbf{e}_1M_1M_2\cdots M_n\mathbf{e}_1^T=0.$ For this purpose, define for $0\le a<k<b,$ $$P_k(a,b,+):=P(\text{starting from } k, \{X_n\} \text{ exits } [a+1,b-1] \text{ from above}). $$
Then one follows from Markov properties that $$P_k(a,b,+)=\frac{\sum_{j=a+1}^k\mathbf{e}_1M_jM_{j+1}\cdots M_{b-1}\mathbf{e}_1^T}{\sum_{j=a+1}^b\mathbf{e}_1M_jM_{j+1}\cdots M_{b-1}\mathbf{e}_1^T},$$ where the empty product  is identity.
Therefore, we have that for $b>k,$ $$P_k(\{X_n\} \text{ hits } k-1 \text{ before it hits } [b,\infty))=1-\frac{\mathbf{e}_1M_kM_{k+1}\cdots M_{b-1}\mathbf{e}_1^T}{1+\sum_{j=k}^{b-1}\mathbf{e}_1M_jM_{j+1}\cdots M_{b-1}\mathbf{e}_1^T}.$$
Since $\lim_{n\rto}\mathbf{e}_1M_1M_2\cdots M_n\mathbf{e}_1^T=0,$ and $M_n,n\ge 0$ are uniformly bounded from below, then
$$P_k(\{X_n\} \text{ hits } k-1 \text{ for some } n\ge0)=1.$$
It follows that $P(T<\infty)=1.$
Hence $\{X_n\}$ is recurrent and so is $\{N_t\}.$

To prove part (b) of the theorem, suppose that $\sum_{n=1}^\infty\frac{1}{\mu_i} \mathbf{e}_1M_1M_2\cdots M_n\mathbf{e}_1^T<\infty.$ Then we have that $\lim_{n\rto}\mathbf{e}_1M_1M_2\cdots M_n\mathbf{e}_1^T=0$ since $\mu_n, n\ge 1$ are uniformly bounded from above. Consequently, it follows from part (a) of the theorem that both $\{X_n\}$ and $\{N_t\}$ are recurrent. Therefore $P(\eta<\infty)=P(T<\infty)=1.$

Let $U_1=\mathbf{e}_1$ and for $i\ge 2$ define
$$U_{i,r}=\#\{0<k< T:X_{k-1}<i,X_k=i+r-1\},\ r=1,...,R.$$
Here and throughout, $``\#\{\ \}"$ denotes the number of elements in set $``\{\ \}".$
Set
$$U_i:=(U_{i,1},U_{i,2},\cdots,U_{i,R}).$$
Then one follows from Hong and Wang \cite{hwa} (see Theorem 1.1 therein)
that $(U_n)_{n\ge1}$ forms an R-type branching process with offspring distribution \begin{eqnarray}\label{up01}
&&P(U_{i+1}=(u_1,...,u_L)\big|U_{i}=\mathbf{e}_1)\no\\
&&\quad\quad=\frac{(u_1+...+u_L)!}{u_1!\cdots
u_L!}\z(\frac{\lambda_i^1}{\mu_i+\sum_{k=1}^R\lambda_i^k}\y)^{u_1}\cdots\z(\frac{\lambda_i^R}{\mu_i+\sum_{k=1}^R\lambda_i^k}\y)^{u_L}\z(\frac{\mu_i}{\mu_i+\sum_{k=1}^R\lambda_i^k}\y),
\end{eqnarray}
and for $2\le l\le L,$
\begin{eqnarray}\label{up02}
&&P\z(U_{i+1}=(u_1,...,1+u_{l-1},...,u_L)\big|U_{i}=\mathbf{e}_l\y)\no\\
 &&\quad\quad=\frac{(u_1+...+u_L)!}{u_1!\cdots
u_L!}\z(\frac{\lambda_i^1}{\mu_i+\sum_{k=1}^R\lambda_i^k}\y)^{u_1}\cdots\z(\frac{\lambda_i^R}{\mu_i+\sum_{k=1}^R\lambda_i^k}\y)^{u_L}\z(\frac{\mu_i}{\mu_i+\sum_{k=1}^R\lambda_i^k}\y).
\end{eqnarray}
It follows from (\ref{up01}) and (\ref{up02}) that
\begin{equation}\label{eu}E(U_i)=\mathbf{e}_1A_1A_2\cdots A_{i-1}\end{equation} where \begin{equation}\label{am}A_i=\left(
    \begin{array}{cccc}
   b_i^1 & ... & b_i^{R-1} & b_i^R \\
   1+b_i^1 & ... &  b_i^{R-1} & b_i^R \\
    \vdots & \ddots & \vdots & \vdots \\
    b_i^1 & ... & 1+ b_i^{R-1} & b_i^R\\
    \end{array}
     \right)\end{equation} with $b_i^r=\frac{\lambda_i^r}{\mu_i},\ r=1,...,R.$
Considering the occupation time of $\{X_n\}$ before $T,$ we have that $\sum_{k=0}^{T-1}1_{X_k=0}=1$ and for $i\ge1,$ \begin{equation}\label{ox}\sum_{k=0}^{T-1}1_{X_k=i}=U_{i,1}+U_{i+1}\mathbf{1}.\end{equation}
Considering the occupation time of $\{N_t\}$ before $\eta,$ we have $\int_{0}^\eta1_{N_t=0}dt=\xi_{0,1}$ and for $i\ge 1,$ \begin{equation}\label{on}\int_{0}^\eta1_{N_t=i}dt=\sum_{k=1}^{U_{i,1}+U_{i+1}\mathbf1}\xi_{i,k},\end{equation} where $\xi_{i,k},i\ge0,k\ge1$ are mutually independent random variables, which are also independent of $U_{i},$ such that $P(\xi_{i,k}>t)=e^{-t(\mu_i+\sum_{r=1}^R\lambda_i^r)}, t\ge 0.$
For the proof of (\ref{ox}) and (\ref{on}), refer to Wang \cite{w}.

By Ward's equation, it follows from (\ref{eu}) and (\ref{on}) that  \begin{equation*}
  \begin{split}
    E(\eta)&=\sum_{n=0}^\infty E\Big(\int_0^\eta 1_{N_t=n}dt\Big)=E\xi_{0,1} +\sum_{n=1}^\infty E(U_{n,1}+U_{n+1}\mathbf1)E\xi_{n,1}\\
    &=\Big(\sum_{r=1}^R\lambda_0^r\Big)^{-1}+\sum_{n=1}^\infty\frac{1}{\mu_n+\sum_{r=1}^R\lambda_n^r}(\mathbf{e}_1A_1A_2\cdots A_{n-1}\mathbf{e}_1^T+\mathbf{e}_1A_1A_2\cdots A_{n}\mathbf1)\\
    &=\Big(\sum_{r=1}^R\lambda_0^r\Big)^{-1}+\sum_{n=1}^\infty\frac{1}{\mu_n}\mathbf{e}_1A_1A_2\cdots A_{n-1}\mathbf1\\
 &=\Big(\sum_{r=1}^R\lambda_0^r\Big)^{-1}+\sum_{n=1}^\infty\frac{1}{\mu_n}\mathbf{e}_1M_1M_2\cdots M_{n-1}\mathbf{e}_1^T<\infty.
  \end{split}
\end{equation*}
Therefore $\{N_t\}$ is positive recurrent. The existence of the limits $\psi_k$ in (\ref{lima}) and (\ref{limb}) follows from Theorem 1.6 of Chapter 5 in Anderson \cite{ads}. In the same theorem, it is also showed that if $(\psi_k)_{k\ge0}$ is a probability distribution, then it is the unique probability distribution such that for all $t>0,$ \begin{equation}\label{iv}\psi_k=\sum_{n=0}^\infty \psi_nP_t(n,k),\ k\ge 0.\end{equation}
On the other hand, in Theorem 3.5.1 of Norris \cite{Nor}, it is shown that $\Big(\frac{\pi_k}{\mu_k+\sum_{r=1}^R\lambda_k^r}\Big)_{k\ge0}$ satisfies (\ref{iv}), where $(\pi_k)_{k\ge0}$ is the stationary distribution of the embedded process $\{X_n\}.$ Next we calculate $(\pi_k)_{k\ge0}.$ We mention that for (1,R) random walk which could stay at its position, $(\pi_k)_{k\ge0}$ was studied in Hong and Zhou \cite{hz}. Since the calculation is not too long and the notations are a bit different from those in \cite{hz}, we repeat the calculation of $(\pi_k)_{k\ge0}$ here.

Note that by (\ref{eu}) and (\ref{ox}),  $E(\sum_{k=0}^{T-1}1_{X_k=0})=1$ and for $n\ge1,$
\begin{equation}
\begin{split}
E\Big(\sum_{k=0}^{T-1}1_{X_k=n}\Big)&=E(U_{n,1}+U_{n+1}\mathbf{1})\\
&= \mathbf{e}_1(A_1A_2\cdots A_{n-1}\mathbf{e}_1^T+A_1A_2\cdots A_n\mathbf1)\\
&=\frac{\mu_n+\sum_{r=1}^R\lambda_n^r}{\mu_n} \mathbf{e}_1A_1A_2\cdots A_{n-1}\mathbf1\\
&=\frac{\mu_n+\sum_{r=1}^R\lambda_n^r}{\mu_n}\mathbf{e}_1M_1M_2\cdots M_{n-1}\mathbf{e}_1^T.
\end{split}
\end{equation}

Also we have
 \begin{equation*}\begin{split}
ET&=1+\sum_{n=1}^\infty E(U_{n,1}+U_{n+1}\mathbf1)\\
  &=1+\sum_{n=1}^\infty\frac{\mu_n+\sum_{r=1}^R\lambda_n^r}{\mu_n} \mathbf{e}_1M_1M_2\cdots M_{n-1}\mathbf{e}_1^T.
\end{split}
\end{equation*}
Let $\pi_0=\frac{1}{1+\sum_{n=1}^\infty\frac{\mu_n+\sum_{r=1}^R\lambda_n^r}{\mu_n} \mathbf{e}_1M_1M_2\cdots M_{n-1}\mathbf{e}_1^T}$ and for $k\ge1$
$$\pi_k=\frac{\frac{\mu_k+\sum_{r=1}^R\lambda_k^r}{\mu_k}\mathbf{e}_1M_1M_2\cdots M_{k-1}\mathbf{e}_1^T}{1+\sum_{n=1}^\infty \frac{\mu_n+\sum_{r=1}^R\lambda_n^r}{\mu_n} \mathbf{e}_1M_1M_2\cdots M_{n-1}\mathbf{e}_1^T}.$$

Then one follows from Durrett \cite{dur} that $(\pi_k)_{k\ge0}$ defines a stationary distribution of $\{X_n\}.$ Set $\nu_k=\pi_k/(\mu_k+\sum_{r=1}^R\lambda_k^r)$ for $k\ge 0.$ Then $(\nu_k)_{k\ge0}$ satisfies (\ref{iv}). Normalizing $\nu_k$ by setting $\psi_{k}=\nu_k/(\sum_{k=0}^\infty\nu_k)$ we get (\ref{lima}) and (\ref{limb}). We conclude that $\psi_k,k\ge0$ is the stationary distribution of $\{N_t\}.$ The theorem is proved. \qed

\noindent{\large{\bf \large Acknowledgements:}} The author would like to thank Dr. Hongyan Sun for her useful discussion when writing the paper.


\end{document}